\documentclass{JHEP3}
\usepackage{amsmath,amsfonts,amssymb,amsthm,epsfig}
\usepackage{graphicx}

\input xy
    \xyoption{matrix}
    \xyoption{arrow}
    \xyoption{curve}

\def\calP{{\mathcal P}}
\def\calG{{\mathcal G}}
\def\calL{{\mathcal L}}

\newcommand{\N}{{\mathbb N}}

\def\PG{\mathrel{\calP\calG}}
\def\bar#1{\overline{#1}}
\def\obj{{\rm obj}}
\def\mor{{\rm mor}}
\def\Aut{{\rm Aut}}
\def\Groups{{\bf Groups}}
\def\at#1#2{\left. {#1} \right|_{#2}}

\def\normal{\triangleleft}
\def\To{\longrightarrow}
\def\Mapsto{\longmapsto}
\def\te#1{\text{#1}}
\def\eqref#1{(\ref{#1})}

\def\f{\frac}

\newtheorem{theorem}{Theorem}

\theoremstyle{definition}

\theoremstyle{remark}

\newtheorem{conjecture}{Conjecture}

\def\proofthm#1{\vskip 0.05 in
        \noindent {\it Proof of Theorem \ref{#1}.} $\ \ $ }

\def\endproof{$\Box$}

\def\be{\begin{equation}}
\def\ee{\end{equation}}
\def\beq{\begin{equation}}
\def\eeq{\end{equation}}

\def\bea{\begin{eqnarray}}
\def\eea{\end{eqnarray}}

\newcommand\beano{\begin{eqnarray*}}
\newcommand\eeano{\end{eqnarray*}}
\newcommand\ben{\begin{enumerate}}
\newcommand\een{\end{enumerate}}

\newenvironment{superbox}{
    \begin{center}
    \begin{tabular}{|c|}\hline\\
    \begin{minipage}{12cm}
}{
    \end{minipage}\\\\
    \hline
    \end{tabular}
    \end{center}
}

\title{Geometric Characterization of Property R}

\author{{ William Gordon Ritter}\\
{\small Harvard University Department of Physics} \\
{\small 17 Oxford St., Cambridge, MA 02138} \\
{\it Email: ritter@fas.harvard.edu}}

\abstract{ Consider pairs of the form $(G, N)$, with $G$ a group
and $N \normal G$, as objects of a category $\PG$. A morphism
$(G_1, N_1) \To (G_2, N_2)$ will be a group homomorphism $f : G_1
\To G_2$ such that $f(N_1) \subset N_2$. We introduce a functor $Q
: \PG \To \Groups$, which provides a geometric definition of
Property R, since it is most naturally visualized by means of a
directed graph. We compute these graphs for a number of finite
groups of small order, and prove a general characterization of the
graphs which occur in this way.}


\preprint{}

\begin{document}

\section{Introduction}

For a group $G$, we denote the derived subgroup by either $G'$ or
$[G,G]$. $G$ is said to have \emph{Property R} if $(G/N)' \ne G/N$
holds for every proper normal subgroup $N \normal G$, i.e.
nontrivial quotients of $G$ are never perfect. Interesting
relationships between different normal subgroups can be modeled by
forming graphs with the normal subgroups as points, and edges
isolating the property of interest. The edge relation which
characterizes Property R is given explicitly in
eq.~\eqref{edge-rel}.

All groups which are either solvable, simple, or are symmetric
groups satisfy Property R. The latter property is important mainly
because it appears as a hypothesis in powerful simplicity theorems
for groups with a Tits system or BN-pair. A well known example is
the Tits Simplicity Theorem \cite{Bou}.

\begin{theorem}[Tits Simplicity, \cite{Bou}] \label{tits}
Let $G$ be a group having a BN-pair. Let $Z$ be the intersection
of the conjugates of $B$, let $U < B$ and let $G_1$ be the
subgroup generated by the conjugates of $U$ in $G$. Assume $U
\normal B = UT$, $U$ has Property R, $G_1 = G_1'$, and $(W,S)$ is
an irreducible Coxeter system \footnote{See \cite{Bou}, IV, Sec.~1
no 9 for definitions}. Then any subgroup $H \subset G$ normalized
by $G_1$ is contained in $Z$  or contains $G_1$. It follows that
$G_1/(G_1 \cap Z)$ is either trivial, or is nonabelian and simple.
\end{theorem}

\section{The Category of Pairs}

Let the objects of category $\PG$ be pairs of the form $(G, N)$
where $G$ is a group and $N \normal G$. A morphism between two
objects $(G_1, N_1) \To (G_2, N_2)$ will be a group homomorphism
$f : G_1 \To G_2$ such that $f(N_1) \subset N_2$. For each such
$f$, we denote the induced homomorphism by $\bar{f} : G_1/N_1 \To
G_2/N_2$. In this brief note we define a functor $Q : \PG \To
\Groups$, and study its fundamental properties. This functor
provides a geometric definition of Property R, and is naturally
visualized by means of a directed graph.

For any group $G$, we use the notation
\[
    G' = [G,G]
\]
for the derived subgroup, generated by elements of the form
$xyx^{-1}y^{-1}$. Consider the natural 1-1 correspondence given by
the Lattice isomorphism theorem:
\[
    \left\{
        \begin{matrix}
            \te{ subgroups of } G  \\
            \te{ containing } N
        \end{matrix}
    \right\}
    \longleftrightarrow
    \big\{
        \te{ subgroups of } G/N
    \big\}
\]
\be \label{barmap}
    H \ \ \Mapsto \ \  \bar{H} = H/N  \leq G/N
\ee
Given an object $(G, N) \in \obj_{\PG}$, the commutator subgroup
$[G/N, G/N]$ is a characteristic subgroup in $G/N$. Let $Q(G, N)$
denote the preimage of $[G/N, G/N]$ under the map \eqref{barmap}.
The diagram is as follows:
\[
\xymatrix{
    G \ar@{-}[d] \ar[r]  &  G/N  \ar@{-}[d]  \\
    Q(G, N)  \ar@{-}[d] \ar@{-}[r]   &  [G/N, G/N]  \ar@{-}[d]  \\
    N \ar[r]  &  \{ e \}
}
\]
Given a group $G$ and a normal subgroup $N \normal G$,  $Q$ associates a
second normal subgroup of $G$ which is not, in general, equal to $N$.

\begin{theorem} \label{functor}
$Q$ is a functor on the category $\PG$.
\end{theorem}

\proofthm{functor} Given a morphism $f \in \mor_{\PG}$, $f : (G_1,
N_1) \To (G_2, N_2)$, we define an associated morphism $Q(f) :
Q(G_1, N_1) \To Q(G_2, N_2)$ by restriction:
\[
    Q(f) = \at{f}{Q(G_1, N_1)} \, .
\]
Let $i \in \{1,2\}$ and suppose $G_i, N_i, Q_i$ are arbitrary
groups satisfying $N_i \normal Q_i \normal G_i$ and $f : (G_1,
N_1) \To (G_2, N_2)$ is a morphism as above. If
\[
    \bar{f} \big( Q_1 / N_1 \big) \subset Q_2 / N_2
\]
then $f(Q_1) \subset Q_2$. (If $\exists\, q_1 \in Q_1$ such that
$f(q_1) \not\in Q_2$, then $\bar{f}(q_1 N_1) = f(q_1)N_2 \not\in
Q_2/N_2$).

It remains to see that
\be \label{containment}
    f\big(Q(G_1, N_1)\big)
        \subset
    Q(G_2, N_2)
\ee
This will finish the proof, since it shows that $Q(f)$ maps into
$Q(G_2, N_2)$. To prove \eqref{containment}, note that
\[
    \bar{f}\big( Q(G_1, N_1) / N_1 \big)
    =
    \bar{f}\big( (G_1 / N_1)' \big)
    =
    \big( \bar{f}(G_1 / N_1) \big)'
    \subset
    (G_2 / N_2)'
    =
    Q(G_2, N_2) / N_2
\]
The conclusion now follows by the argument given above.
\endproof

\section{The Forest Associated to $G$}

Throughout this paper, we use the term \emph{tree} to mean graphs
with no cycles except possibly self-loops. $Q$ gives a natural way
of associating a forest of trees to a group $G$. Let $T_G$ be the
digraph with vertices labelled by the normal subgroups of $G$, and
with edge relation $e$ defined as follows:
\be \label{edge-rel}
    N_1 \, e \, N_2 \iff  N_2 = Q(G, N_1)
\ee
We refer to this forest as the \emph{Property R graph} of $G$,
since $G$ has Property R if and only if the vertex in $T_G$
corresponding to the whole group $G$ is isolated. An example of
such a graph for a group of order 192 is the following:

    \begin{center}
    \vskip 0.1in
\epsfig{file=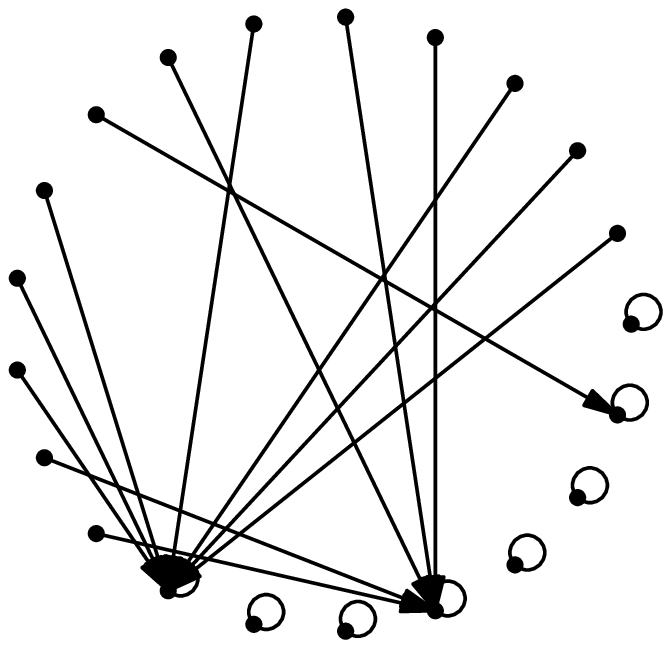,width=1in}
    \vskip 0.1in
    \end{center}

Since the first homology $H_1(G)$ is the same as the
abelianization $G / [G,G]$, our edge relation implies that if
there is an edge $a \to b$, then $H_1(G/a) \cong G/b$. We now
state our main results, and provide proofs where necessary.

\begin{theorem} \label{main}
The graph defined by \eqref{edge-rel} has the following properties:
\end{theorem}
\begin{enumerate}
\item {\bf (Inclusion in the full lattice)}
Let $\calL$ be the full lattice of normal subgroups of $G$, with
edge relation defined by inclusions. Since $N_1 \subseteq Q(G,
N_1)$, it follows that $T_G$ can be viewed as a subgraph $\calL$
with the same vertices and fewer edges.

\item {\bf (Self-loops)}
Self-loops correspond to abelian quotients, in the following sense:
\[
    N = Q(G,N)
    \iff
    G/N \te{ is abelian }
    \iff
    \te{ vertex } N \te{ has a self-loop}
\]
Since $T_G$ has no loops other than self-loops, we can write more
simply in this context
\[
    \te{abelian quotient } = \te{ loop}
\]
A quotient $G/N$ of $G$ is abelian if and only if $N$ includes
$G'$, so the Property R graph provides a quick way of visualizing
the normal subgroups which are between $G$ and $G'$ in the
lattice.

\item {\bf (Property R)}
The graph associated to $Q$ gives a geometric definition of
Property R. $G$ has Property R if and only if the vertex in $T_G$
corresponding to the whole group $G$ is isolated, because a normal
subgroup $N$ with $(G/N)' = G/N$ generates a directed edge of the
form $N \To G$.

\item {\bf (Automorphism group action)}
There is a natural action of the automorphism group $\Aut(G)$ on
the tree $T_G$, i.e. if $N_1 \To N_2$ is an edge in $T_G$, then
$f(N_1) \To f(N_2)$ is also an edge in $T_G$, for any $f \in
\Aut(G)$.

\emph{Proof:} $\ $ Since $N_2/N_1 \subseteq \big( G/N_1 \big)'$,
we may represent any coset $n_2 N_1 \in N_2/N_1$ as $g' N_1$,
where $g' \in G'$. Then $f(n_1 N_2) = f(g') f(N_1) \in \big(
G/f(N_1) \big)'$. We have used the property that $G'$ is
characteristic, and hence preserved under all automorphisms. The
reverse inclusion follows similarly.

\item {\bf (Edge stabilizers)}
An arbitrary edge of $T_G$ may be written $N \To Q(G,N)$. The
stabilizer of this edge under the action of the group $\Aut(G)$ on
the forest $T_G$ is the same as the stabilizer of $N$ under the
action of $\Aut(G)$ on subgroups of $G$.
\end{enumerate}

\section{Characterization of Property R Graphs}

\newcommand{\st}{{\cal S}}
\newcommand{\selfloop}{\xymatrix{ \circ \ar@(ur,ul)[] }}

The $n$-\emph{star graph}, denoted $S_n$, is a tree on $n+1$ nodes
with one node having vertex degree $n$ and the others having
degree 1.

    \begin{center}
    \vskip 0.1in
\epsfig{file=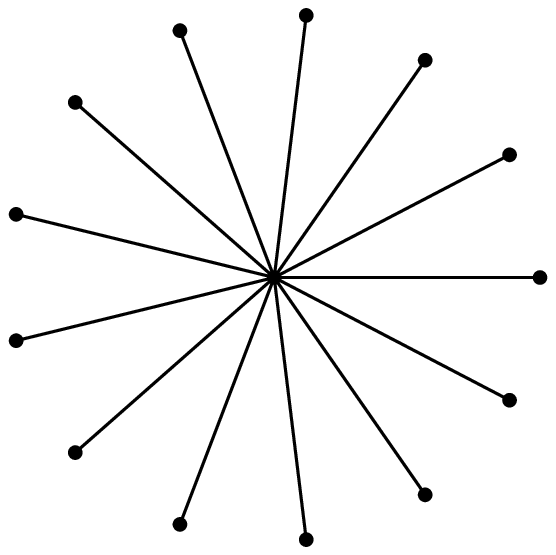,width=1in}
    \vskip 0.1in
    \end{center}

We also define $\st_n$ to be the star graph $S_n$, modified by the
addition of a self-loop at the center. It is easy to see that
Property R graphs, modulo self-loops, are disjoint unions of
stars.

\begin{theorem} \label{lemma:stars}
The Property R graph of any group $G$ is a union of disjoint
copies of $\st_j$ for various values of $j \in \N$, and isolated
points. If we define $\st_{-1}$ to be an isolated point, then the
isomorphism type of the graph is encoded in one list of numbers,
$\{j : \st_j \subset T, j \geq -1\}$, counted with multiplicity.
\end{theorem}

\proofthm{lemma:stars} Suppose to the contrary that there is a
nontrivial configuration of the form $a \To b \To c$ where $a,b,c$
are distinct proper normal subgroups. The arrow $a \To b$ means
that $b/a \cong (G/a)'$. Then we have
\[
    G/b \cong \frac{G/a}{b/a} \cong \f{G/a}{(G/a)'} = H_1(G/a)
\]
which is abelian. Therefore $(G/b)' \cong \{e\}$ and so any arrow
from $b$ must be a self-loop. In particular $b = c$. $\Box$

\vskip 0.1 in

Our computations lead us to conjecture the converse of Theorem
\ref{lemma:stars}:

\begin{conjecture}
For any configuration ${\cal G}$ of stars $\st_j$ and isolated
points as in Lemma \ref{lemma:stars}, there exist some group $G$
such that $T_G = {\cal G}$.
\end{conjecture}

\section{ Computation and Specific Examples }

For small groups, the Property R graphs $T_G$ studied in the
previous section can be computed explicitly and visualized through
simple and elegant computer algebra methods. We include a table of
such computations for finite groups of small order. These
calculations were performed using \emph{GAP} for the group theory,
and \emph{Mathematica} for the combinatorics and graph
visualization. The GAP program analyzes the lattice of normal
subgroups, computes the necessary quotients and derived groups,
and generates graphs in a format which is acceptable input for
\emph{Combinatorica}, a standard package for Mathematica since
version 4.

The heart of the group theory program is the following four lines:
    \vskip 0.1 in

\begin{superbox} \tt \small
 NS := NormalSubgroups(G);

 for i in [ 1 .. Length(NS) ] do

 $~~~$phi := NaturalHomomorphismByNormalSubgroupNC(G, NS[i]);

 $~~~$j := Position(NS, PreImage(phi, DerivedSubgroup(Image(phi))));
 \rm
\end{superbox}

    \vskip 0.1 in
The GAP program at this point writes instructions for
\emph{Combinatorica} to place a directed edge from vertex $i$ to
vertex $j$ in the graph it is building.

To denote some specific examples, we introduce the notation that
if $G, H$ are two graphs, then the product $GH$ denotes the
disjoint union. In particular, $G^n = \coprod_{i=1}^n G$.
Eq.~\eqref{GL} shows the result of some computations for general
linear groups of low order.
\be \label{GL}
    \begin{matrix}
        GL(2,7) & (\st_0 \st_2)^2 \\
        GL(2,13) & (\st_0 \st_1 \st_2)^3 \\
        GL(2,17) & \st_0 \st_2 (\st_1)^3 \\
        GL(2,19) & (\st_0 \st_2)^3
    \end{matrix}
\hskip 1.3 in
    \begin{matrix}
        GL(3,2) & \st_1  \\
        GL(3,3) & (\st_1)^2 \\
        GL(3,4) & \st_0 \st_2 \\
        GL(3,5) & (\st_1)^3
    \end{matrix}
\ee
In addition, we have the following isomorphisms between graphs,
\[
    T_{GL(2,23)} = T_{GL(3,7)} = T_{GL(2,7)} = (\st_0 \st_2)^2
\]
The Property R graphs of simple groups, such as $PSL$ over a
finite field, for obvious reasons all take the form $\bullet\!\!
\To \!\!\bullet\!\!\bigcirc$, which is denoted $\st_1$. It is not
an interesting property for simple groups.

Figure \ref{gl2} and Figure \ref{gl3} give the Property R diagrams
for $2 \times 2$ and $3 \times 3$ matrix groups over finite fields
$GF(q^n)$ for the smallest interesting values of $q$ and $n$.
\footnote{The author wishes to thank the Harvard Math department
for some processor time on a powerful Sun computer, which was used
for the calculations in Figure \ref{gl2} and Figure \ref{gl3}.}

\FIGURE{
    \centering
    \includegraphics{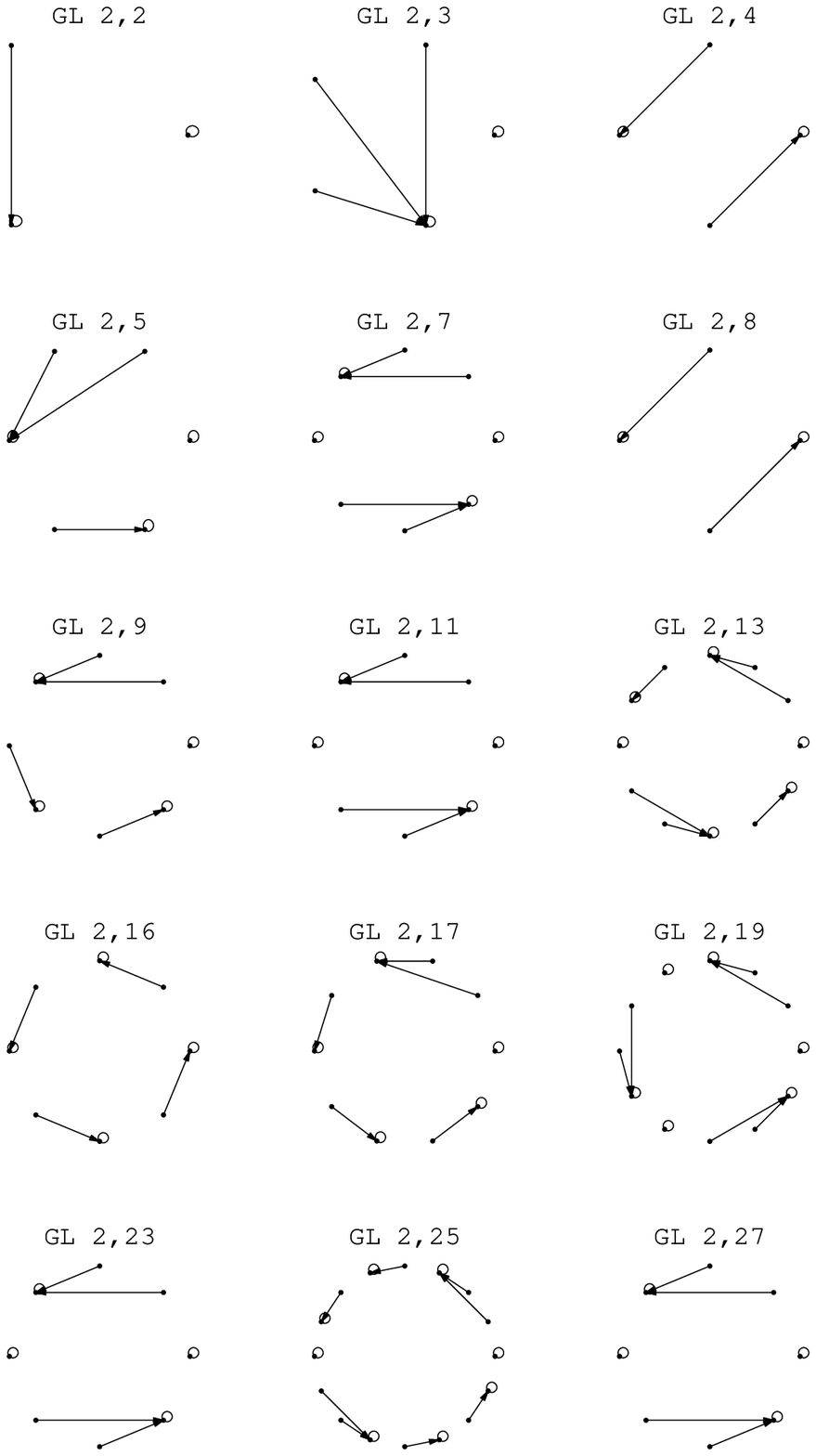}
    \caption{$GL(2,q^n)$}
    \label{gl2}
}

\FIGURE{
    \centering
    \includegraphics{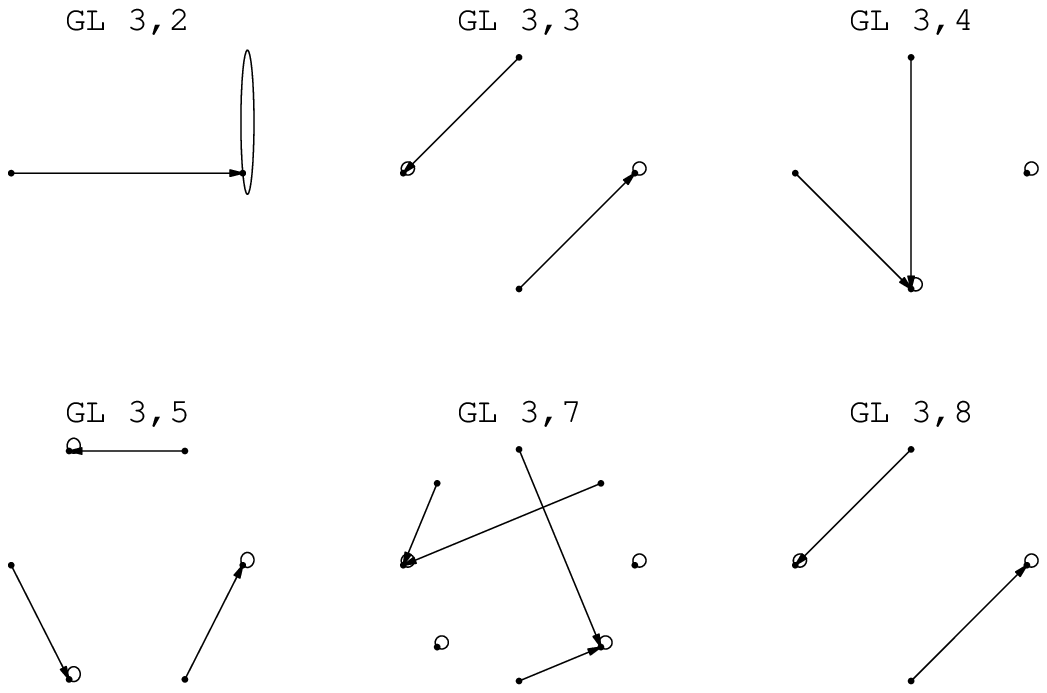}
    \caption{$GL(3,q^n)$}
    \label{gl3}
}

\end{document}